\documentclass[10pt]{amsart}
\usepackage{amsthm,amsmath,amssymb,tikz}
\usetikzlibrary{arrows}
\usepackage{lscape}

\title[Near-central enumerative methods]{Character-theoretic techniques for near-central enumerative problems}
\author{D. M. Jackson$^1$ and C. A. Sloss$^2$}

\theoremstyle{plain}
\newtheorem{thm}{Theorem}[section] 
\newtheorem{lemma}[thm]{Lemma} 
 
\newtheorem{cor}[thm]{Corollary} 
\newtheorem{prob}[thm]{Problem}

\theoremstyle{definition} 
\newtheorem{defn}[thm]{Definition}

\newcommand{\bbC}{\mathbb{C}}
\newcommand{\bul}{\bullet}
\newcommand{\Ga}{\Gamma}
\newcommand{\frH}{\mathfrak{H}}
\newcommand{\la}{\lambda}
\newcommand{\La}{\Lambda}
\newcommand{\frS}{\mathfrak{S}}

\newcommand{\SYT}{\mathrm{SYT}}

\renewcommand{\includegraphics}{}



\thanks{
${\hspace{-1ex}}^1$ Department of Combinatorics and Optimization, University of Waterloo, Waterloo, Ontario, Canada. Partially supported by an NSERC Discovery Grant. \texttt{dmjackson@math.uwaterloo.ca}}

\thanks{
${\hspace{-1ex}}^2$ Department of Combinatorics and Optimization, University of Waterloo, Waterloo, Ontario, Canada. Partially supported by an NSERC Postgraduate Scholarship. \texttt{csloss@theorem.ca}}

\begin{document}

\parskip =  2pt 

\maketitle

\begin{abstract}
The centre of the symmetric group algebra $\bbC[\frS_n]$ has been used successfully for studying important problems in enumerative combinatorics. These include maps in orientable surfaces and ramified covers of the sphere by curves of genus $g$, for example. However, the combinatorics of some equally important $\frS_n$-factorization problems forces $k$ elements in $\{1,\ldots,n\}$ to be distinguished. Examples of such problems include the star factorization problem, for which $k=1,$ and the enumeration of $2$-cell embeddings of dipoles with two distinguished edges~\cite{VisentinWieler:2007} associated with Berenstein-Maldacena-Nastase operators in Yang-Mills theory~\cite{ConstableFreedmanHeadrick:2002}, for which $k=2.$  Although distinguishing these elements obstructs the use of central methods, these problems may be encoded algebraically in the centralizer of $\bbC[\frS_n]$ with respect to the subgroup $\frS_{n-k}.$  We develop methods for studying these problems for $k=1,$ and demonstrate their efficacy on the star factorization problem.  In a subsequent paper~\cite{JacksonSloss:2011}, we consider a special case of the the above dipole problem by means of these techniques. 
\end{abstract}

\section{Introduction}

\subsection{Background}
The character theory of the symmetric group $\mathfrak{S}_n$ is a powerful tool in enumerative combinatorics. It has been used, for example, to study the enumeration of maps in orientable surfaces (see \cite{JacksonVisentin:1990} and \cite{JacksonVisentin:2001}) and the enumeration of ramified covers of the sphere by a curve of genus $g$ (see \cite{GouldenJackson:1997} and \cite{GouldenJacksonVakil:2005}). An overview of the key elements in the character-based approach to enumerative problems may be found in \cite{Jackson:1987}. The approach relies on having an encoding of an enumerative problem in the centre $Z(n)$ of the complex group algebra $\mathbb{C}[\mathfrak{S}_n]$. 
Let $\kappa(\pi)$ denote the cycle type of an element $\pi\in \mathfrak{S}_n$. Thus, $\kappa(\pi)$ is a partition of the integer $n$, and we write $\kappa(\pi) \vdash n$. The set
$\left\{K_\la \colon \la\vdash n\right\}$
is a basis of $Z(n)$, where $K_\la$ is the formal sum of all elements in
$\mathcal{C}_{\la} := \{\pi\in \frS_n : \kappa(\pi) =\la\}$, 
the conjugacy class of $\frS_n$ naturally indexed by $\la$. 
Hence, the character-based approach is limited to problems for which the solution depends only on the conjugacy class of the permutations involved. 

Problems which cannot be encoded in $Z(n)$ will be referred to as \textit{non-central} and they constitute an important emerging class of problems.  The following algebra plays a role in non-central enumerative problems analogous to that of $Z(n)$  for central ones. Let $\frH$ be a subgroup of 
$\frS_n$. The \textit{centralizer of} $\bbC[\frS_n]$ \textit{with respect to} $\frH$ is the set 
$$Z_{\frH}(n) := \{g\in \bbC[\frS_n] : \sigma g \sigma^{-1} = g \text{ for all } \sigma\in \frH\}.$$
As a notational shorthand, if $\frH = \frS_{\{1,\ldots, k\}}$ for some $k\leq n$, $Z_\frH(n)$ is denoted by $Z_{n-k}(n)$. Clearly $Z_{\frH}(n)$ is an algebra, and  $Z_0(n)$ is the centre of $\bbC[\frS_n]$. 
The centralizers of $\bbC[\frS_n]$ play a key role in an alternative derivation of the irreducible representations of the symmetric group due to Okounkov and Vershik \cite{OkounkovVershik:1996}, who give a general expression for a standard basis for centralizers of the form $Z_k(n)$. The centralizer relevant to the present purpose is $Z_1(n)$,  which has \textit{standard basis} 
$$\left\{   K_{\la,i} \colon \la\vdash n, i\in \la \right\}$$
where
$K_{\la,i} := \sum_{\pi \in \mathcal{C}_{\la,i}} \pi $ and
$\mathcal{C}_{\la,i} := \{\sigma \in \mathcal{C}_{\la} : n \text{ is on a cycle of length } i\}.$ We write $a\in \la$ to indicate that $a$ is a part of $\la$. The orbit-stabilizer theorem immediately gives
$$|\mathcal{C}_{\la,i}| = \frac{(n-1)! \,i\, m_i(\la)}{\prod_i i^{m_i(\la)} m_i(\la)!},$$
where $m_i(\lambda)$ is the multiplicity of $i$ as a part of $\lambda$. A combinatorial problem which may be encoded as an element of $Z_1(n)$ is referred to as a \textit{near-central problem}. 

Let $\mathrm{Inv}$ denote the linear operator on $\bbC[\frS_n]$ defined by $\mathrm{Inv}(\pi) = \pi^{-1}$. The fact that $\mathrm{Inv}(K_{\la,i}) = K_{\la,i}$ gives an elementary proof that $Z_1(n)$ is commutative: for $G,H\in Z_1(n)$, we have
$GH = \mathrm{Inv}(GH) = \mathrm{Inv}(H) \mathrm{Inv}(G) = HG$.
In contrast, when $k\geq 2$, $Z_k(n)$ is non-commutative.

\subsection{The star factorization problem and non-centrality}
The non-central combinatorial problem which is the subject of this paper, and which will be used for developing a general approach to non-central enumerative problems through the use of $Z_1(n)$, is the following. It concerns ordered factorizations $(\tau_1,\ldots,\tau_r)$ of a permutation $\pi$ into transpositions $\tau_i$, where each $\tau_i$ is of the form $(j,n)$ for some $1\leq j\leq n-1$. A transposition of this kind is a \textit{star transposition}. 
\begin{prob}[Star Factorization Problem]
\label{problem:star-factorizations}
Let $\pi\in \frS_n$ and let $r\geq 1$. Determine the number of sequences $(\tau_1,\ldots,\tau_r) \in \frS_n^r$ such that
$\tau_1\tau_2\cdots\tau_r= \pi,$
where each $\tau_i$ is a star transposition.
\end{prob}
This paper demonstrates how algebraic methods, analogous to the character-based methods used to study central problems, can be applied to the star factorization problem.

\subsubsection{The transitive star factorization problem}
A closely related problem is the \textit{Transitive Star Factorization Problem}, in which there is the additional condition:
\emph{the group $\left\langle\tau_1,\ldots,\tau_r\right\rangle$ acts transitively on $\{1,\ldots,n\}$.} This was introduced by Pak \cite{Pak:1998}, who solved the problem for the special case in which $n=km+1$ for some $k$ and $m$, $\pi(n)=n$, $\pi$ is of cycle type $(k^m,1)$, and in which the number of factors is minimal. In this case, transitivity is forced and therefore was not cited as a condition.
 Irving and Rattan \cite{IrvingRattan:2006} solved the minimal case in which the permutation $\pi$ is arbitrary. The transitive version of the problem was solved in full generality by Goulden and Jackson \cite{GouldenJackson:2007}, who gave a solution which, unexpectedly, depends only on the cycle type of $\pi$. In other words the transitive star factorization problem is central.

\subsubsection{Product of Jucys-Murphy elements}
It is easily seen that Problem~\ref{problem:star-factorizations} is equivalent to determining the $r^{\text{th}}$ power of the \textit{Jucys-Murphy element} 
$$J_n := \sum_{1\leq i<n} (i,n).$$
An indirect approach to Problem~\ref{problem:star-factorizations} is through the centrality of symmetric polynomials evaluated at Jucys-Murphy elements~\cite{Jucys:1974} together with an expression for power sums of Jucys-Murphy elements found by Lascoux and Thibon \cite{LascouxThibon:2004}, and the observation that 
$$J_n^r = \sum_{2\leq k\leq n} J_k^r - \sum_{2\leq k \leq n-1} J_k^r.$$
Combinatorially, this amounts to regarding the problem as the \textit{difference} of two central problems.

Our perspective is to approach the problem directly as a near-central problem, and then use an algebra that respects near-centrality, namely $Z_1(n).$  Determining $J_n^r$ is then regarded as a special case of determining products of standard basis elements of $Z_1(n)$.  The techniques employed in this paper are therefore applicable to any near-central problem which has an encoding in terms of standard basis elements.

\subsection{The $(p,q,n)$-dipole problem}
Another significant non-central problem is the problem of enumerating dipoles embedded in an orientable surface with two distinguished edges, in which the separation of the ends of these two edges is specified by parameters $p$ and $q$. These dipoles are the summation indices of two-point functions of Berenstein-Maldacena-Nastase operators, and arise in the study of duality between Yang-Mills theory and string theory conducted by Constable \textit{et al.}~\cite{ConstableFreedmanHeadrick:2002}. They obtained asymptotic results for the torus and double torus. Subsequently, Visentin and Wieler \cite{VisentinWieler:2007} determined exact formulas for the torus and double torus.  In general, the distinguishing of two edges forces the $(p,q,n)$-dipole problem to lie in $Z_2(n)$. However, a subsequent paper \cite{JacksonSloss:2011} shows that when $q=n-1$, the $(p,q,n)$-dipole problem is near-central. The techniques developed in the present paper are then applied to give a solution to the $(p,n-1,n)$-dipole problem on all orientable surfaces. 

\subsection{Organization of the paper}
Section \ref{sec:definitions} lists the definitions and results from the enumerative theory of $Z(n)$ which are used in this paper. Section \ref{sec:centralizers} introduces the algebraic context for studying the star factorization problem, namely, the centralizer of $\bbC[\frS_n]$ with respect to the subgroup $\frS_{n-1}$. The connection coefficients of this algebra are given in terms of generalized characters, which were introduced by Strahov \cite{Strahov:2007}. Section \ref{sec:evaluating-GC} shows how explicit formulae for certain generalized characters may be obtained, generalizing a technique of Diaconis and Greene \cite{DiaconisGreene:1989}, and establishes relationships between sums of generalized characters and irreducible characters of the symmetric group. Section \ref{sec:SF-formulae} gives an expression for the solution to the star factorization problem in terms of generalized characters, and gives some special cases where more explicit solutions may be obtained.

\section{Notation, definitions and background results}
\label{sec:definitions}
In this section we review the standard results which we make use of in this paper (see also \cite{JamesKerber:1981, Macdonald:1995}). Let $\mathbb{P}$ denote the set of positive integers. Let $\la =(\la_1, \la_2, \ldots)$ be a partition of $n$.  Then each $\la_i$ is called a \textit{part} of $\la$. If $i$ is a part which is repeated $a_i$ times, it is convenient to write it more succinctly as $i^{a_i}$. Recall that the \textit{multiplicity} $m_i(\la)$ of $i>0$ in $\la$ is the number of times $i$ appears in $\la$. Let $m(\la)$ denote the total number of parts of $\la$. It will be convenient to have notation for the following partitions constructed from $\la$:  \\
--- if $i\in\la,$ then $\la \setminus i$ is the partition obtained by reducing $m_i(\la)$ by one; \\[2pt]
--- for any $i,$  $\la\cup i$ is the partition obtained by increasing $m_i(\la)$ by one; \\[2pt]
--- $i_-(\la)$ is the partition obtained by replacing a part $i$ by $i-1.$

Let $\frS_n$ denote the symmetric group acting on $\{1,\ldots,n\}$, and let $\mathbb{C}[\mathfrak{S}_n]$ denote its complex group algebra. The right-to-left convention for multiplication in $\mathfrak{S}_n$ will be used, \textit{i.e.} $\pi_1\pi_2$ is the permutation obtained by applying $\pi_2$ followed by $\pi_1$. Let $m(\pi)$ denote the number of cycles of $\pi$. Let $\chi^{\la}_{\mu}$ denote the irreducible character of $\mathfrak{S}_n$ indexed by the partition $\la$, evaluated at a permutation of cycle type $\mu$. Let $d_{\la}$ denote the degree of the irreducible representation indexed by $\la$ and let
\[
X^{\la} :=  \frac{d_{\la}}{n!} \sum_{\mu\vdash n} \chi^{\la}_{\mu} K_{\mu}.
\]
Then $\left\{X_\la \colon \la\vdash n\right\}$ is a basis of orthogonal idempotents of $Z(n).$ 

Let $\mathcal{F}_{\la}$ be the Ferrers diagram of shape $\la\vdash n$, and let $\mathrm{SYT}_{\la}$ be the set of all standard Young tableaux of shape $\la$ on $\{1,\ldots,n\}$ (English convention). Let $T^\bul$ be the tableau obtained by deleting the cell containing $n$ from $T$.
The \textit{content} $c_T(i)$ of $i$ in a standard Young tableau $T$ is defined to be $k-j$ where $i$ is in the cell in row $j$ and column $k$.  The \textit{content vector} of $T$ is
$\mathbf{c}_T := (c_T(1), c_T(2),\ldots, c_T(n))$. In a slight abuse of notation, $\mathbf{c}_{\lambda}$ denotes the multiset of contents of any tableau of shape $\lambda\vdash n$. The \textit{content polynomial} corresponding to the partition $\lambda$ is  
\[
c_{\lambda}(t) = \prod_{1\leq i\leq n} (t + c_i(T))
\]
for any $T\in \mathrm{SYT}_{\la}$.  

Given $T\in \mathrm{SYT}_{\la}$, let $e_T\in \mathbb{C}[\mathfrak{S}_n]$ denote Young's semi-normal unit corresponding to $T$. This paper makes use of the following facts about the semi-normal units.
\begin{lemma}
\label{lem:seminormal-facts}
Let $\la\vdash n$. Then the following results hold. 

\noindent (1) 
The degree of the irreducible representation indexed by $\la$ is
$d_{\la} = |\mathrm{SYT}_{\la}|.$ \\
\noindent (2)
If $T \in \mathrm{SYT}_{\la}$  then the coefficient of the identity in $e_T$ is $d_{\la}/n!$. \\
\noindent (3) 
Let $T,S\in \mathrm{SYT}_{\la}.$  Then $e_Te_S = \delta_{T,S}e_T$. \\
\noindent (4)
Let $\la \vdash n$. Then
$X^{\la} = \sum_{T\in \mathrm{SYT}_{\la}} e_T.$ \\
\noindent (5)
If $T$ is a tableau with $n$ boxes, then $e_{T^\bul} = \sum_{S : S^\bul=T^\bul} e_S$. \\
\noindent (6)
Let $T\in \mathrm{SYT}_{\la}$ and let $k\leq n$. Then
$J_k \, e_T = c_T(k)\, e_T.$
\end{lemma}
The quantity $d_{\lambda}$ may be determined in general using the hook-length formula, although for the partitions $\lambda$ arising in this paper, $d_{\lambda}$ may be determined by elementary combinatorial arguments. For example,
\[
d_{(n-k,1^k)} = \binom{n-1}{k},
\]
since a tableau of shape $(n-k,1^k)$ is determined by the choice of $k$ symbols to put in the first column, from among the symbols $\{2,\ldots,n\}$. 

A polynomial $f \in \mathbb{C}[x_1,\ldots,x_n]$ is said to be \textit{symmetric} if, for any $\pi\in \mathfrak{S}_n$, $f(x_1,x_2,\ldots,x_n) = f(x_{\pi(1)}, x_{\pi(2)},\ldots, x_{\pi(n)})$. The notation $\Lambda[x_1,\ldots,x_n]$ denotes the ring of symmetric polynomials in the indeterminates $x_1,\ldots,x_n$. Two classes of symmetric polynomials are used in this paper. The \textit{elementary symmetric polynomials} are defined by $e_k := [t^k]E(t)$, where $E(t):=\prod_{1\leq i\leq n} (1+tx_i)$. The \textit{power sum symmetric polynomials} are defined by $p_k := \sum_{1\leq i\leq n} x_i^k$.

\section{Algebraic methodology for near-central problems}
\label{sec:methodology}

\subsection{Connection coefficients for $Z_1(n)$}
\label{sec:centralizers}

A set of orthogonal idempotents for $Z_1(n)$ may be defined by restricting the summation in the expression for $X^{\la}$ given in part (4) of Lemma \ref{lem:seminormal-facts}. Let $\mathrm{SYT}_{\la,i}$ denote the set of standard Young tableaux of shape $\la$ in which $n$ appears at the end of a row of length $i$. Let $c_{\la,i}$ denote the content of the box containing $n$ in any tableau in $\mathrm{SYT}_{\la,i}$. The quantity $c_{\la,i}$ depends only on $\la$ and $i$, and is given by
\[
c_{\la,i} = i - \sum_{k \geq i} m_k(\la).
\]

\begin{defn}[$Z_1$-idempotents]
\label{defn:Z1-idempotents}
Let $\la\vdash n$ and let $i$ be a part of $\la$. The $Z_1$\textit{-idempotents} are the elements in $\mathbb{C}[\mathfrak{S}_n]$ given by
\[
\Ga^{\la,i} := \sum_{T \in \mathrm{SYT}_{\la,i}} e_T.
\]
\end{defn}
Several facts about $\Ga^{\la,i}$ are immediate from this definition and from properties of the semi-normal units, and are thus stated below without proof.
\begin{lemma}
Let $\la,\mu\vdash n$. Let $i\in \la$ and $j\in \mu$. Then the following statements are true.
\begin{enumerate}
\item
$\sum_{k\in \la} \Ga^{\la,k} = X^{\la}.$ 
\item  $\Ga^{\la,i} \Ga^{\mu,j} = \delta_{\la,\mu} \delta_{i,j} \, \Ga^{\la,i}. $ 
\item
The set $\{\Ga^{\la,i} \colon \la\vdash n, i\in \la \}$ is linearly independent. 
\item
$\vert \{\Ga^{\la,i} \colon \la\vdash n, i\in \la \} \vert =  
\vert \{K_{\la,i} \colon \la\vdash n, i\in \la \} \vert.$ 
\end{enumerate}
\label{lem:Z1-idempotent-properties}
\end{lemma}

It is not immediately obvious from the definition of the $Z_1$-idempotents that they do lie in $Z_1(n)$. This may be proven as follows.
\begin{lemma}
\label{lemma:Gamma-in-Z1}
Let $\la\vdash n$ and $i$ be a part of $\la$. Then
$
\Ga^{\la,i} = X^{\la}X^{i_-(\la)},
$
where the product is taken in $\mathbb{C}[\mathfrak{S}_n]$ and every element in the support of $X^{i_-(\la)}$ is regarded as having the element $n$ as a fixed point. Consequently, $\Ga^{\la,i}\in Z_1(n)$, and $\{\Ga^{\la,i}\}_{\la\vdash n, i\in\la}$ is a basis for $Z_1(n)$.
\end{lemma}
\begin{proof}
By part 4 of Lemma \ref{lem:seminormal-facts},
\[
X^{\la}X^{i_-(\la)} = \left(\sum_{T \in \SYT_{\la}} e_T  \right)\left( \sum_{S \in \SYT _{i_-(\la)}} e_S \right).
\]
For a fixed $S\in \SYT_{i_-(\la)}$, apply parts (3) and (5) of Lemma \ref{lem:seminormal-facts} to obtain
\begin{align*}
\sum_{T\in \SYT_{\la}} e_T e_S &= \sum_{T\in \SYT_{\la}} \sum_{\substack{S_0\in \SYT_{\la}, \\ S_0^{\bul}=S}} e_T e_{S_0} 
= \sum_{\substack{T\in \SYT_{\la}, \\ T^\bul= S}} e_T. 
\end{align*}
Summing over all $S \in \SYT_{i_-(\la)}$ gives
\[
X^{\la}X^{i_-(\la)} = \sum_{S \in \SYT_{i_-(\la)}} \, \sum_{\substack{T\in \SYT_{\la}, \\ T^\bul= S}} e_T.
\]
Since the tableaux $T\in \SYT_{\la}$ with the property that $T^\bul\in \SYT_{i_-(\la)}$ are those in which the symbol $n$ appears at the end of a row of length $i$, then
\[
X^{\la}X^{i_-(\la)}  = \sum_{T\in \SYT_{\la,i}} e_T = \Ga^{\la,i}.
\]
\end{proof}

Since $\Ga^{\la,i}\in Z_1(n)$, then the coefficients of $\Ga^{\la,i}$ in the standard basis are well-defined:
\begin{defn}
\label{defn:GC-definition-1}
Let $\la,\mu\vdash n$ and let $i$ and $j$ be parts of $\la$ and $\mu$, respectively. Define
\[
\gamma^{\la,i}_{\mu,j} := \frac{n!}{d_{\la}} [K_{\mu,j}] \Ga^{\la,i}.
\]
\end{defn}

The coefficients $\gamma^{\la,i}_{\mu,j}$ were first studied by Strahov \cite{Strahov:2007}, who named them \textit{generalized characters} of the symmetric group. (Although the term ``generalized character'' is also used in the literature to refer to any integral linear combination of characters, this is not the sense of the term used here.) Strahov defined generalized characters as the zonal spherical functions of the Gel'fand pair $(\mathfrak{S}_n\times \mathfrak{S}_{n-1}, \mathrm{diag}(\mathfrak{S}_{n-1}))$, where the diagonal subgroup $\mathrm{diag}(G)$ of $G \times G$ is defined by
$\mathrm{diag}(G) = \{(g,g) : g\in G\}.$ 
Strahov provides the following expression for generalized characters in terms of irreducible characters of the symmetric group for any $\pi\in\mathcal{C}_{\mu,j}$:
\[
\gamma^{\la,i}_{\mu,j} = \frac{d_{i_-(\la)}}{(n-1)!} \sum_{\sigma\in \mathfrak{S}_{n-1}} \chi^{\la}(\pi\sigma^{-1})\chi^{i_-(\la)}(\sigma).
\]
This expression may also be obtained routinely by extracting coefficients from the equation $\Ga^{\la,i} = X^{\la}X^{i_-(\la)}$. Thus, Lemma \ref{lemma:Gamma-in-Z1} is a proof that Definition \ref{defn:GC-definition-1} is equivalent to Strahov's definition of generalized characters. Consequently, the coefficients $\gamma^{\la,i}_{\mu,j}$ satsify the properties of zonal spherical functions. (A list of these properties may be found in Chapter VII of Macdonald \cite{Macdonald:1995}.) Of particular relevance to the present task is the fact that these coefficients are orthogonal with respect to the standard inner product on $\mathbb{C}[\mathfrak{S}_n]$. This is the content of the following result.
\begin{cor}
\label{cor:GC-orthogonality}
Let $\la,\mu \vdash n$ and let $i\in\la$, $j\in\mu$. Then
\[
\frac{1}{n!} \sum_{\substack{\rho \vdash n,\\ k\in \rho}} |\mathcal{C}_{\rho,k}| \gamma^{\la,i}_{\rho,k}\gamma^{\mu,j}_{\rho,k} =\frac{d_{i_-(\la)}}{d_{\la}} \delta_{\la,\mu}\delta_{i,j}
\]
\end{cor}
Since the coefficients of $\Ga^{\la,i}$ in the standard basis are orthogonal, inverting the change-of-basis transformation is routine, giving the following.
\begin{lemma}
\label{lemma:standard-to-idempotent}
Let $\la \vdash n$ and let $i$ be a part of $\la$. Then
\[
K_{\la,i} = \sum_{\mu\vdash n, j\in \mu} \frac{|\mathcal{C}_{\la,i}|}{d_{j_-(\mu)}} \gamma^{\mu,j}_{\la,i} \Ga^{\mu,j}.
\]
\end{lemma}

The preceding results allow the connection coefficients for $Z_1(n)$ to be expressed in terms of generalized characters. Let $\la,\mu\vdash n$ and let $i\in\la$ and $j\in\mu.$  The product $K_{\la,i}K_{\mu,j}$ may be written in the $Z_1$-idempotent basis as
\[
K_{\la,i}K_{\mu,j} = |\mathcal{C}_{\la,i}||\mathcal{C}_{\mu,j}| \sum_{\la'\vdash n, i'\in \la} \, \sum_{\mu'\vdash n, j'\in \mu'} \frac{\gamma^{\la',i'}_{\la,i}}{d_{i'_-(\la')}} \frac{\gamma^{\mu',j'}_{\mu,j}}{d_{j'_-(\mu')}} \Ga^{\la',i'} \Ga^{\mu',j'}.
\]
By orthogonal idempotency of the $Z_1$-idempotents, this may be written as 
\[
K_{\la,i} K_{\mu,j} = |\mathcal{C}_{\la,i}||\mathcal{C}_{\mu,j}|  \sum_{\rho \vdash n, \ell \in \rho} \frac{\gamma^{\rho,\ell}_{\la,i} \gamma^{\rho,\ell}_{\mu,j}}{d_{\ell_-(\rho)}^2} \Ga^{\rho,\ell}.
\]
Extracting the coefficient of $K_{\nu,k}$ yields the following.
\begin{thm}
\label{lemma:Z1-connection-coefficients}
Let $\la,\mu,\nu\vdash n$ and let $i,j$ and $k$ be parts of $\la,\mu$ and $\nu$, respectively. Let the coefficients $c_{\la,i,\mu,j}^{\nu,k}$ be defined by
$
K_{\la,i}K_{\mu,j} = \sum_{\nu\vdash n, k\in \nu} c_{\la,i,\mu,j}^{\nu,k} K_{\nu,k}.
$
Then
\[
c_{\la,i,\mu,j}^{\nu,k} = \frac{ |\mathcal{C}_{\la,i}||\mathcal{C}_{\mu,j}| }{n!} \sum_{\rho \vdash n, \ell \in \rho} \frac{\gamma^{\rho,\ell}_{\la,i} \gamma^{\rho,\ell}_{\mu,j}\gamma^{\rho,\ell}_{\nu,k}}{d_{\ell_-(\rho)}} \frac{d_{\rho}}{d_{\ell_-(\rho)}}.
\]
\end{thm}

This result has the following generalization to the case involving an arbitrary number of factors, which may be proven in a similar manner.
\begin{thm}
\label{thm:Z1-r-factors}
For $1\leq k\leq r$, let $\lambda^{(k)}\vdash n$ and let $i_k\in \lambda^{(k)}$. Let $\mu\vdash n$ and $j\in \mu$. Then
\[
[K_{\mu,j}] \prod_{1\leq k\leq r} K_{\lambda^{(k)}, i_k} = \frac{1}{n!} \sum_{\substack{\rho\vdash n \\ \ell\in \rho}} \frac{\gamma^{\rho,\ell}_{\mu,j}d_{\rho}}{d_{\ell_-(\rho)}^r} \prod_{1\leq k\leq r} |\mathcal{C}_{\lambda^{(k)},i_k}| \gamma^{\rho,\ell}_{\lambda^{(k)},i_k}.
\]
\end{thm}

\subsection{Evaluation of Generalized Characters}
\label{sec:evaluating-GC}

The method of finding explicit formulae for generalized characters presented here generalizes a method due to Diaconis and Greene \cite{DiaconisGreene:1989} for evaluating irreducible characters of the symmetric group. Their method relies on the following fact. 
\begin{lemma}
\label{lemma:diaconis-greene}
Let $\la\vdash n$ and let $f \in \La[x_2,\ldots,x_n]$ be such that $K_{\la} = f(J_2,\ldots,J_n)$. Then
$$
\chi^{\mu}_{\la} = \frac{d_{\mu}}{|\mathcal{C}_{\la}|} f(\mathbf{c}(\mu)).
$$
\end{lemma}
The following generalization of this result holds for generalized characters. Let $\La^{(1)}[x_2,\ldots,x_n]$ denote the ring polynomials that are invariant under permutations of $x_2,\ldots,x_{n-1}$. Such a polynomial will be called an \textit{almost symmetric} polynomial, and may be regarded as a polynomial in $x_n$ whose coefficients are symmetric polynomials in the variables $x_2,\ldots, x_{n-1}$. The generalization of the Diaconis-Greene approach relies on the following. 
\begin{lemma}
\label{lemma:genchar-from-asf}
Let $\mu\vdash n$ and let $j$ be a part of $\mu$. Let $f\in \La^{(1)}[x_2,\ldots,x_n]$ be such that $f(J_2,\ldots,J_n) = \sum_{\substack{\la\vdash n \\ i \in \la}} a_{\la,i} K_{\la,i}$. Then $\Ga^{\mu,j}$ is an eigenvector of $f(J_2,\ldots,J_n)$ with eigenvalue $f(\mathbf{c}_{j_-(\mu)},c_{\mu,j})$, and 
\[
\sum_{\substack{\la\vdash n, \\ i \in \la}} a_{\la,i} \frac{|\mathcal{C}_{\la,i}|}{d_{j_-(\mu)}} \gamma^{\mu,j}_{\la,i} = f(\mathbf{c}_{j_-(\mu)}, c_{\mu,j}).
\]
\end{lemma}
\begin{proof}
First, observe that
$1 = K_{(1^n),1} = \sum_{\mu\vdash n,  j\in \mu} \Ga^{\mu,j}.$
Indeed, this follows from Lemma \ref{lemma:standard-to-idempotent} and part 2 of Lemma \ref{lem:seminormal-facts}, since
\[
\gamma^{\mu,j}_{(1^n),1} = \frac{n!}{d_{\mu}} [K_{(1^n),1}] \Ga^{\mu,j} =  |\SYT_{\mu,j}| = d_{j_-(\mu)}. 
\]
Thus,
\begin{align*}
\sum_{\substack{\la\vdash n, \\ i \in \la}} a_{\la,i} K_{\la,i} &= f(J_2,\ldots,J_n) K_{(1^n),1} 
= \sum_{\substack{\mu\vdash n \\ j\in \mu}}  f(J_2,\ldots,J_n)  \Ga^{\mu,j}.
\end{align*}
By the definition of $\Ga^{\mu,j}$,
\begin{align*}
f(J_2,\ldots,J_n)  \Ga^{\mu,j} &= \sum_{T\in \SYT_{\mu,j}} f(J_2,\ldots,J_n)  e_T \\
&=  \sum_{T\in \SYT_{\mu,j}} f(c_T(2),c_T(3),\ldots, c_T(n)) e_T \\
&= \sum_{T\in \SYT_{\mu,j}} f( \mathbf{c}_{j_-(\mu)}, c_{\mu,j}) e_T.
\end{align*}
The quantity  $f( \mathbf{c}_{j_-(\mu)}, c_{\mu,j})$ is well-defined since $f$ is symmetric in $x_2,\ldots,x_{n-1}$. 
Thus,
\[
f(J_2,\ldots,J_n)  \Ga^{\mu,j}  =  f( \mathbf{c}_{j_-(\mu)}, c_{\mu,j}) \Ga^{\mu,j},
\]
and
\[
\sum_{\substack{\la\vdash n, \\ i \in \la}} a_{\la,i} K_{\la,i} =  \sum_{\substack{\mu\vdash n \\ j\in \mu}}  f( \mathbf{c}_{j_-(\mu)}, c_{\mu,j}) \Ga^{\mu,j}.
\]
On the other hand, by Lemma \ref{lemma:standard-to-idempotent}, the standard basis for $Z_1(n)$ may also be expressed in the generalized character basis as follows:
\[
\sum_{\substack{\la\vdash n, \\ i \in \la}} a_{\la,i}K_{\la,i} = \sum_{\substack{\la\vdash n \\ i \in \la}} a_{\la,i}\sum_{\substack{\mu\vdash n \\ j\in \mu}} \frac{|\mathcal{C}_{\la,i}|}{d_{j_-(\mu)}} \gamma^{\mu,j}_{\la,i} \Ga^{\mu,j}.
\]
Comparing coefficients in the $Z_1$-idempotent basis gives the result.
\end{proof}
An immediate corollary of this result is the following expression for generalized characters.
\begin{cor}
\label{cor:genchar-from-asf}
Let $\la,\mu\vdash n$ and let $i$ and $j$ be parts of $\la$ and $\mu$, respectively. Let $f\in \Lambda^{(1)}[x_2,\ldots,x_n]$ be such that $f(J_2,\ldots,J_n) = K_{\la,i}$. Then 
 the generalized characters are given by the formula
\[
\gamma^{\mu,j}_{\la,i} = \frac{d_{j_-(\mu)}}{|\mathcal{C}_{\la,i}|}  f( \mathbf{c}_{j_-(\mu)}, c_{\mu,j}).
\]
\end{cor}
This result may be used to express the generalized character $\gamma^{\la,i}_{\mu,j}$ as an evaluation of an almost symmetric polynomial, provided there is an explicit expression for $K_{\mu,j}$ as an almost symmetric polynomial in Jucys-Murphy elements. (The existence of such a polynomial is guaranteed by a result of Olshanski  \cite{Olshanski:1988}, namely, that $\Lambda^{(1)}[J_2,\ldots,J_n] = Z_1(n)$.) Some examples of generalized characters which can be evaluated using this method are given in Table \ref{table:GC-evaluations} of Appendix \ref{sec:appendix-A}, using the expressions for basis elements of $Z_1(n)$ appearing in Table \ref{table:Z1-JM-basis}. 

There are two important relationships between sums of generalized characters and irreducible characters of the symmetric group. The first follows directly from Definitions \ref{defn:Z1-idempotents} and \ref{defn:GC-definition-1}:
\begin{lemma}
\label{lemma:genchar-superscript-sum}
Let $\la,\mu\vdash n$.  
Then $\chi^{\mu}_{\la} =  \sum_{j\in \mu} \gamma^{\mu,j}_{\la,i}$
for any $i\in \la$.
\end{lemma}
Lemma \ref{lemma:genchar-from-asf} may also be used to prove an analogous relationship in which the sum is taken over subscripts, instead of superscripts. 
\begin{lemma}
\label{lemma:genchar-subscript-sums}
Let $\la, \mu\vdash n$ and let $j$ be a part of $\mu$. Then
\[
\chi^{\mu}_{\la} = \frac{d_{\mu}}{|\mathcal{C}_{\la}| d_{j_-(\mu)}} \sum_{i\in\la} |\mathcal{C}_{\la,i}| \gamma^{\mu,j}_{\la,i}.
\]
\end{lemma}
\begin{proof}
Let $f_{\la,i}\in \Lambda^{(1)}[x_2,\ldots,x_n]$ be such that $f_{\la,i}(J_2,\ldots,J_n) = K_{\la,i}$, and let $f_{\la}\in \Lambda[x_2,\ldots,x_n]$ be a symmetric polynomial in $x_2,\ldots,x_n$ such that $f_{\la}(J_2,\ldots,J_n) = K_{\la}$. Thus,
\[
\sum_{i\in \la} f_{\la,i}(J_2,\ldots,J_n) = f_{\la}(J_2,\ldots,J_n).
\]
For any $\mu\vdash n$ and $j\in \mu$, let $T\in \SYT_{\mu,j}$. Then
\begin{align*}
f_{\la}(\mathbf{c}_{\mu}) e_T &= f_{\la}(J_2,\ldots, J_n) e_T 
= \sum_{i\in \la} f_{\la,i}(J_2,\ldots,J_n) e_T  \\
&= \sum_{i\in \la} f_{\la,i}(\mathbf{c}_{j_-(\mu)}, c_{\mu,j}) e_T,
\end{align*}
so
$f_{\la}(\mathbf{c}_{\mu})  = \sum_{i\in \la} f_{\la,i}(\mathbf{c}_{j_-(\mu)}, c_{\mu,j}).$
Applying Lemmas \ref{lemma:diaconis-greene} and \ref{lemma:genchar-from-asf} yields the result.
\end{proof}
These results may be used to obtain some generalized character formulae which do not appear in Table \ref{table:GC-evaluations}. For example,
\begin{cor}
\label{cor:genchar-1-fixed-not-n}
Let $\mu \vdash n$ and let $j$ be a part of $\mu$. Then
\[
\gamma^{\mu,j}_{(n-1,1),n-1} = 
\begin{cases}
1 & \text{~if } \mu=(n), j=n, \\
(-1)^n & \text{~if } \mu = (1^n), j=1, \\
\frac{(-1)^{k+1}}{n-1} & \text{~if } \mu = (n-k,1^k), 1\leq k\leq n-2, j=n-k, \\
\frac{(-1)^k}{n-1} & \text{~if } \mu = (n-k,1^k), 1\leq k\leq n-2, j=1, \\
\frac{(-1)^k}{k(n-k-2)} & \text{~if } \mu = (n-k-1,2,1^{k-1}), j=2, \\
0 & \text{~otherwise.}
\end{cases}
\]
\end{cor}
\begin{proof}
Rearranging the equation in Lemma \ref{lemma:genchar-subscript-sums} gives
\[
\gamma^{\mu,j}_{(n-1,1),n-1} = \frac{1}{n-1} \left( \frac{n d_{j_-(\mu)}}{d_{\mu}} \chi^{\mu}_{(n-1,1)} - \gamma^{\mu,j}_{(n-1,1),1}\right).
\]
Substituting the known value for and $\gamma^{\mu,j}_{(n-1,1),1}$ along with the fact that
\[
\chi^{\mu}_{(n-1,1)} = 
\begin{cases}
1 & \text{~if } \mu = (n), \\
(-1)^n & \text{~if } \mu = (1^n), \\
(-1)^k & \text{~if } \mu = (n-k-1,2,1^{k-1}), \\
0 & \text{~otherwise}
\end{cases}
\]
gives the result. (The formula for $\chi^{\mu}_{(n-1,1)}$ is a known result, but it can also be obtained using Lemma \ref{lemma:genchar-superscript-sum} and Table \ref{table:GC-evaluations}.)
\end{proof}

Summing generalized characters over all partitions with a fixed number of parts gives a more explicit result.
\begin{lemma}
\label{lem:gen-cha-weighted-sum}
Let $\rho\vdash n$ and let $\ell$ be a part of $\rho$. Then
\[
\sum_{\substack{\la\vdash n, \\ m(\la)= m, i\in \la}} \frac{|\mathcal{C}_{\la,i}|}{d_{\ell_-(\rho)}} \gamma^{\rho,\ell}_{\la,i} =  e_{n-m}(\mathbf{c}_{\rho}) = [t^{m}] c_{\rho}(t).
\]
\end{lemma}
\begin{proof}
It is routine to show that the elementary symmetric polynomial $e_{n-m}(x_2,\ldots,x_n)$ evaluated at the Jucys-Murphy elements is given by 
\[
e_{n-m}(J_2,\ldots,J_n) = \sum_{\substack{\la\vdash n \\ m(\la)= m}} K_{\la} =  \sum_{\substack{\la\vdash n \\ m(\la)= m, i\in \la}} K_{\la,i}.
\]
Thus, by Lemma \ref{lemma:genchar-from-asf},
\[
\sum_{\substack{\la\vdash n \\ m(\la)= m,i\in \la}}\frac{|\mathcal{C}_{\la,i}|}{d_{\ell_-(\rho)}} \gamma^{\rho,\ell}_{\la,i} = e_{n-m}(\mathbf{c}_{\rho}).
\]
The result then follows by using the generating series for elementary symmetric functions. 
\end{proof}

\section{Applications to the Star Factorization Problem}
\label{sec:SF-formulae}

In this section, the methods of Section \ref{sec:methodology} are applied to the study of the star factorization problem. The most general form of the solution to the star factorization problem, in terms of generalized characters, is given in Theorem \ref{lemma:JM-powers-as-gen-char}. Knowledge of explicit expressions for generalized characters and of algebraic relationships between generalized characters allows us to identify special cases of the star factorization problem which have more structure than the general problem. A new formula is given for the number of star factorizations of a permutation $\pi\in \mathcal{C}_{(n-1,1),n-1}$ in Corollary \ref{cor:JM-factorizations-full-cycle}. The number of star factorization of all $\pi$ in a given conjugacy class is given in terms of \textit{ordinary} irreducible characters in Theorem \ref{thm:NTSF-character-based}. Finally, Theorem \ref{thm:NTSF-cycle-sum} establishes a relationship between the set of all star factorizations of a permutation with $r$ cycles and content polynomials.

The material in Section \ref{sec:methodology} provides the following two methods for determining $J_n^r$.
\begin{enumerate}
\item
Use Theorem \ref{thm:Z1-r-factors}, with $\lambda^{(k)} = (2,1^{n-2})$ and $i_k = 2$ for $1\leq k\leq r$, along with the expression
\[
\gamma^{\mu,j}_{(2,1^{n-2}),2} = (n-1)^{-1} d_{j_-(\mu)}c_{\mu,j}
\] 
from Table \ref{table:GC-evaluations}.
\item
Use Lemma \ref{lemma:genchar-from-asf}, with $f(x_2,\ldots,x_n) = x_n^r$. Expressing the identity element in the $Z_1$-idempotent basis gives
\[
J_n^r = f(J_2,\ldots,J_n) \sum_{\substack{\mu\vdash n \\ j\in \mu}} \Gamma^{\mu,j} = \sum_{\substack{\mu\vdash n \\ j\in \mu}} c_{\mu,j}^r \Gamma^{\mu,j}.
\]
\end{enumerate}
Using either method leads to the following result.

\begin{thm}
\label{lemma:JM-powers-as-gen-char}
Let $\la \vdash n$ and let $i$ be a part of $\la$. For $\pi\in \mathcal{C}_{\la,i}$, the number of factorizations of $\pi$ into $r$ star transpositions is given by
\[
[K_{\la,i}]J_n^r = \sum_{\substack{\mu\vdash n \\ j\in \mu}} \frac{d_{\mu}}{n!} \gamma^{\mu,j}_{\la,i} c^r_{\mu,j}.
\]
\end{thm}

Theorem \ref{lemma:JM-powers-as-gen-char} can be used to give more explicit expressions for the coefficients of $J_n^r$ in cases when the generalized characters arising in the expression for $J_n^r$ can be evaluated. For example,  Corollary \ref{cor:genchar-1-fixed-not-n} gives a simple expression for $\gamma_{(n-1,1),n-1}^{\mu,j}$. Thus, Theorem \ref{lemma:JM-powers-as-gen-char} gives the following expressions for the number of factorizations of a permutation in $\mathcal{C}_{(n-1,1),n-1}$ into star transpositions. 
\begin{cor}
\label{cor:JM-factorizations-full-cycle}
Let $r\geq 1$. The number of factorizations of a permutation $\pi \in \mathcal{C}_{(n-1,1),n-1}$ into $r$ star transpositions is given by
\begin{equation}
\label{eqn:JM-factorizations-C}
[K_{(n-1,1),n-1}]J_n^r = \frac{r!}{n!(n-1)} [x^r] \left( 2n\cosh((n-1)x) - 2^n \sinh\left(\frac{(n-1)x}{2}\right) \sinh\left(\frac{x}{2}\right)^{n-1}\right).
\end{equation} 
when $n$ is even, and 
\begin{equation}
\label{eqn:JM-factorizations-D}
[K_{(n-1,1),n-1}]J_n^r = \frac{r!}{n!(n-1)} [x^r] \left( 2n\sinh((n-1)x) - 2^n \sinh\left(\frac{(n-1)x}{2}\right) \sinh\left(\frac{x}{2}\right)^{n-1}\right)
\end{equation}
when $n$ is odd.
\end{cor}
\begin{proof}
By Table \ref{table:GC-evaluations}, the only marked partitions $(\mu, j)$ which make a non-zero contribution to the formula of Theorem \ref{lemma:JM-powers-as-gen-char} are those of the form $((n-k,1^k),n-k)$, $((n-k,1^k),1)$ and $((n-k-1,2,1^{k-1}),2)$. For marked partitions of the form $((n-k-1,2,1^{k-1}),2)$, the quantity $c_{\mu,j}$ is zero, so when $r\geq 1$, these partitions do not contribute to the sum. There are four cases giving non-zero contributions to $[K_{(n-1,1),n-1}]J_n^r$. The pair $(n),n$ contributes
\[
\frac{1}{n!} (n-1)^r.
\]
The pair $(1^n),1$ contributes
\[
\frac{1}{n!} (-1)^n (1-n)^r.
\]
The partitions of the form $(n-k,1^k),n-k$ contribute
\[
\sum_{1\leq k\leq n-2} \binom{n-1}{k} \frac{(-1)^{k+1}}{n!(n-1)} (n-k-1)^r,
\]
and partitions of the form $(n-k,1^k),1$ contribute 
\[
\sum_{1\leq k\leq n-2} \binom{n-1}{k} \frac{(-1)^k}{n!(n-1)} (-k)^r.
\]
Using the fact that $j^k = k![x^k]e^{jx}$, the sum over these four cases simplifies to
\[
[K_{(n-1,1),n-1}] J_n^r = \frac{r!}{n!(n-1)}[x^r]\left(ne^{(n-1)x} + (-1)^n ne^{-(n-1)x} + (1 - e^{(n-1)x})(1-e^{-x})^{n-1} \right).
\]
(Constant terms may be disregarded when extracting the coefficient of $x^r$, since $r\geq 1$.) The result now follows.

\end{proof}

Table \ref{table:GC-evaluations} gives simple expressions for $\gamma_{(n),n}^{\mu,j}$ and $\gamma_{(n-1,1),1}^{\mu,j}$, which may be used to obtain the following expressions in a manner similar to the proof of Corollary \ref{cor:JM-factorizations-full-cycle}.
The number of factorizations of a full cycle of length $n$ into $r$ star transpositions is given by
\begin{equation}
\label{eqn:JM-factorizations-A}
[K_{(n),n}] J_n^r = \frac{2^n (r+1)!}{n!(n-1)}[x^{r+1}] \sinh\left( \frac{(n-1)x}{2}\right) \sinh\left( \frac{x}{2}\right)^{n-1}, 
\end{equation}
and the number of factorizations of a permutation $\pi\in \mathcal{C}_{(n-1,1),1}$ into $r$ star transpositions is given by
\begin{equation}
\label{eqn:JM-factorizations-B}
[K_{(n-1,1),1}]J_n^r = \frac{2^n r!}{n!} [x^r]\sinh\left(\frac{(n-1)x}{2}\right) \sinh\left(\frac{x}{2}\right)^{n-1} .
\end{equation}
Since star factorizations of a elements of $\mathcal{C}_{(n),n}$ and $\mathcal{C}_{(n-1,1),1}$ are necessarily transitive, the formula given in Equations (\ref{eqn:JM-factorizations-A}) and  (\ref{eqn:JM-factorizations-B}) coincide with the formula given by Goulden and Jackson for the same problem. Equation (\ref{eqn:JM-factorizations-C}) is notable for corresponding to a case in which not all the factorizations which are enumerated by this formula are transitive, and which cannot be obtained from the formula of Goulden and Jackson. 

Lemma \ref{lemma:genchar-subscript-sums} may be used to give an expression for the number of star factorizations of all permutations in a given conjugacy class as a linear combination of irreducible characters. This surprising result is an example of an expression involving ordinary irreducible characters of the symmetric group which has been obtained as a result of working in the algebra $Z_1(n)$. It is not known how to obtain this result by working only in the centre of $\mathbb{C}[\mathfrak{S}_n]$. 
\begin{thm}
\label{thm:NTSF-character-based}
Let $\la\vdash n$, and $r\geq 1$. The number of factorizations $(\tau_1,\ldots,\tau_r)$ such that each $\tau_i$ is a star transposition and $\prod_{1\leq i\leq r}\tau_i \in \mathcal{C}_{\la}$ is given by
\[
\frac{|\mathcal{C}_{\la}|}{n!} \sum_{\mu\vdash n}\left( \sum_{j\in \mu} d_{j_-(\mu)} c_{\mu,j}^r \right) \chi_{\la}^{\mu} .
\]
\end{thm}
\begin{proof}
The number of sequences satisfying the stated properties is
\[
\sum_{i \in \la} |\mathcal{C}_{\la,i}| [K_{\la,i}] J_n^r.
\]
By Theorem \ref{lemma:JM-powers-as-gen-char}, after interchanging the order of summation, this expression is equal to
\[
\sum_{\substack{\mu\vdash n, \\ j\in \mu}} \frac{d_{\mu}}{n!} c^r_{\mu,j}\sum_{i \in \la} |\mathcal{C}_{\la,i}| \gamma^{\mu,j}_{\la,i}.
\]
By Lemma \ref{lemma:genchar-subscript-sums}, for any $j\in \mu$,
\[
\sum_{i\in \la} |\mathcal{C}_{\la,i}| \gamma_{\la,i}^{\mu,j} = \frac{|\mathcal{C}_{\la}| d_{j_-(\mu)}}{d_{\mu}} \chi_{\la}^{\mu},
\]
from which the result follows.
\end{proof}

As a further specialization, Lemma \ref{lem:gen-cha-weighted-sum} suggests that the formula for products of Jucys-Murphy elements has a particularly elegant form, in terms of content polynomials, when we are only concerned with the number of cycles of permutations appearing in the product, as opposed to their cycle type. In other words, Lemma \ref{lem:gen-cha-weighted-sum} gives the solution to this problem as a linear combination of content polynomials. 
\begin{thm}
\label{thm:NTSF-cycle-sum}
The number of ordered factorizations $\pi=\tau_1\tau_2\cdots\tau_r$ such that $\pi$ has $k$ cycles and $\tau_i = (j, n)$ for some $1\leq j \leq n-1$ is given by
\[
\sum_{\substack{\mu\vdash n, \\ j\in \mu}} \frac{d_{\mu} d_{j_-(\mu)} c_{\mu,j}^r}{n!} e_{n-k}(\mathbf{c}_{\mu}) = [t^k] \sum_{\substack{\mu\vdash n,\\ j\in \mu}}  \frac{d_{\mu} d_{j_-(\mu)} c_{\mu,j}^r}{n!} c_{\mu}(t).
\]
\end{thm}
\begin{proof}
This number is given by
\begin{align*}
\sum_{\substack{\la \vdash n \\ m(\la)=k}} \sum_{i\in \la} |\mathcal{C}_{\la,i}|[K_{\la,i}] J_n^r &= \sum_{\substack{\la \vdash n \\ m(\la)=k}} \sum_{i\in \la} |\mathcal{C}_{\la,i}| \sum_{\substack{\mu\vdash n \\ j\in \mu}} \frac{d_{\mu}}{n!} \gamma^{\mu,j}_{\la,i} c^r_{\mu,j} \\
&= \sum_{\substack{\mu\vdash n \\ j\in \mu}} \frac{d_{\mu} d_{j_-(\mu)} c_{\mu,j}^r}{n!} \left(\sum_{\substack{\la \vdash n \\ m(\la)=k}} \sum_{i\in \la} \frac{|\mathcal{C}_{\la,i}|}{d_{j_-(\mu)}} \gamma^{\mu,j}_{\la,i} \right),
\end{align*}
from which the result follows by an application of Lemma \ref{lem:gen-cha-weighted-sum}. 
\end{proof}

\appendix

\section{Generalized Character Expressions Obtained Using the Generalized Diaconis-Greene Method}
\label{sec:appendix-A}

In the following tables, for $\la\vdash n$, $\sigma(\la)$ is the sum of contents of a tableau of shape $\la$, and $\sigma^{(2)}$ is the sum of squares of contents of a tableau of shape $\la$.

\begin{table}[h]
\[
\begin{array}{c|c}
\la,i & K_{\la,i} = \\
\hline
(2,1^{n-2}),2 & J_n \\
(2,1^{n-2}),1 & p_1(J_2,\ldots,J_{n-1})  \\
(3,1^{n-3}),3 & J_n^2 - (n-1)K_{(1^n),1} \\
(2,2,1^{n-4}),2 & p_1(J_2,\ldots, J_{n-1})J_n - J_n^2 + (n-1) K_{(1^n),1}  \\
(3,1^{n-3}),1 & p_2(J_2,\ldots,J_{n-1}) - \binom{n-1}{2} K_{(1^n),1} \\
(2,2,1^{n-4}),1 & \frac12(p_1(J_2,\ldots,J_{n-1})^2 - 3p_2(J_2,\ldots,J_{n-1})) + \binom{n-1}{2} K_{(1^n),1}  \\
(n),n & e_{n-1}(J_2,\ldots,J_n) \\
(n-1,1),1 & e_{n-2}(J_2,\ldots,J_{n-1}) 
\end{array}
\]
\caption{Expressions for standard basis elements of $Z_1(n)$ as almost-symmetric polynomials in Jucys-Murphy elements}
\label{table:Z1-JM-basis}
\end{table}

\begin{table}[h]
\[
\begin{array}{c|c}
\la,i &  \gamma^{\mu,j}_{\la,i} = \\
\hline
(2,1^{n-2}),2 & (n-1)^{-1} c_{\mu,j} d_{j_-(\mu)} \\
(2,1^{n-2}),1  & \binom{n-1}{2}^{-1} \sigma(j_-(\mu)) d_{j_-(\mu)} \\
(3,1^{n-3}),3  & \frac12 \binom{n-1}{2}^{-1} (c_{\mu,j}^2-n+1) d_{j_-(\mu)}\\
(2,2,1^{n-4}),2 & (n-1)^{-1}\binom{n-1}{2}^{-1} (\sigma(j_-(\mu))c_{\mu,j} - c_{\mu,j}^2 + n-1) d_{j_-(\mu)} \\
(3,1^{n-3}),1  &  \frac12 \binom{n-1}{3}^{-1} \left(\sigma^{(2)}(j_-(\mu)) - \binom{n-1}{2}\right)d_{j_-(\mu)}\\
(2,2,1^{n-4}),1 & \frac16 \binom{n-1}{4}^{-1} (\sigma(j_-(\mu))^2 - 3\sigma^{(2)}(j_-(\mu)) + (n-1)(n-2)) d_{j_-(\mu)} \\
(n),n &  \begin{cases}
(-1)^{k}\frac{n-k-1}{n-1} & \text{~if } \mu = (n-k,1^{k}), j=n-k; \\
(-1)^{k}\frac{k}{n-1} &\text{~if } \mu = (n-k,1^{k}), j=1; \\
0 & \text{~otherwise.}
\end{cases}\\
(n-1,1),1 &\begin{cases}
(-1)^{k} & \text{~if } \mu = (n-k-1,2,1^{k-1}) \text{ and } j=2; \\
(-1)^{k} & \text{~if } \mu = (n-k,1^{k}) \text{ and } j=n-k; \\
(-1)^{k-1} & \text{~if } \mu = (n-k,1^{k}) \text{ and } j=1; \\
0 & \text{~otherwise.}
\end{cases}
\end{array}
\]
\caption{Explicit formulae for various generalized characters. In this table, $\sigma(\la)$ denotes the sum of contents of any tableau of shape $\la$, and $\sigma^{(2)}(\la)$ denotes the sum of squares of any tableau of shape $\la$.}
\label{table:GC-evaluations}
\end{table}

\bibliographystyle{amsplain}
\bibliography{near-central}

\end{document}